\documentclass[12pt]{article}
\usepackage{fullpage}
\usepackage{amsmath, amsfonts, url}
\newcommand{\lgt}[1]{\log_{10}^{}\,#1\,}
\newcommand{\lgb}[1]{\log_{\,b}^{}#1\,}
\renewcommand{\d}{{\rm d}}
\newcommand{\e}{{\rm e}}
\newcommand{\Blaw}{Benford's law }

\renewcommand{\Pr}[1]{{\mathbb P}\,(#1)}

\newcommand{\PX}[1]{P_X^{}(#1)}
\newcommand{\PS}[1]{P_S^{}(#1)}
\newcommand{\PY}[1]{P_Y^{}(#1)}
\newcommand{\sinc}[1]{{\rm sinc}\,(#1)}
\begin{document}
\begin{center}
   {\Large {\bf Continuous Distributions on $(0,\,\infty)$ Giving} \\
   {\bf Benford's Law Exactly}}\\[2mm]
   Kazufumi OZAWA\,(Professor emeritus at Akita Prefectural
   University, JAPAN)\\
   E-mail: kazufumi.ozawa@gmail.com\\
\end{center}
\begin{quote}
   \centerline{{\bf Abstract}} 
   Benford's law is a famous law in statistics which states that the
   leading digits of random variables in diverse data sets appear not
   uniformly from 1 to 9; the probability that $d\,(d=1,\ldots,9)$
   appears as a leading digit is given by $\log_{10}(1+1/d)$.  This
   paper shows the existence of a random variable with a smooth
   probability density on $(0,\infty)$ whose leading digit
   distribution follows Benford's law exactly.  To construct such a
   distribution the error theory of the trapezoidal rule is used.
\end{quote}
\section{Introduction}
Benford's law, also called the logarithmic law, is a statistical law
which asserts that the leading digits from diverse and different data
sets are not uniformly distributed, as might be expected, but follow
some peculiar distribution described by a logarithmic function.
According to the law, in the sets which obey the law, the number 1
appears as the leading digit about 30\% of the time, while 9 appears
less than 5\% of the time \cite{benford,newcomb}.  \Blaw arises
ubiquitously in the fields of chemistry, physics, geology, astronomy,
engineering, econometrics, etc, and has a wide range of applications
including image processing, network analysis and fraud detection in
accounting data \cite{online, kossovsky}.

It is shown that the distribution of the leading digits of the product
$\prod_{i=1}^n X_i$ converges to Benford's law, as $n\to \infty$,
under the assumption that $X_i$ are independent, identically
distributed and not purely atomic random variables \cite{berger_hill}.
It is also shown that if the distribution of a random variable $X>0$
is the log-normal then the distribution of the leading digits of $X$
approaches to the law as $\sigma^2\to \infty$, where $\sigma^2$ is the
variance of $Y=\ln X$ \cite{berger_hill}.  

The purpose of this note is to construct a continuous distribution on
$(0,\,\infty)$ such that the leading digits of the variable with this
distribution follow \Blaw exactly, without any limiting operations
such as the above.

The rest of the paper is organized as follows.  Section 2 is a brief
introduction of the law.  Section 3 investigates the property of
continuous distributions on $(0,\infty)$ that satisfy the law exactly.
In Section 4, such distribution is proposed by using an error theory
of the trapezoidal rule in numerical integrations.  Section 5 is the
conclusion.

\section{\Blaw}
Let $X$ be a real random variable and its scientific notation or decimal
floating-point representation be
\begin{equation}
    X=\pm S(X)\cdot 10^M,\qquad 1\leq S(X) <10,\quad M\in {\mathbb Z},
	   \label{Eq:scientific} 
\end{equation}
and the decimal notation of the significand $S(X)$ be
\begin{equation*}
	S(X)= d_0.d_1d_2\cdots, 
\end{equation*}
where
\[
    d_0^{}\in \{1,2,\ldots,9\},\qquad d_i^{}\in
    \{0,1,2,\ldots,9\},~~i=1,2,\ldots. 
\]
\Blaw \cite{benford, newcomb} states that if the data $X$ are distributed
over several orders of magnitude, then the leading digit $d_0^{}$ will
approximately follow the probability distribution
\begin{equation}
   \Pr{\,d_0=d\,}=\lgt{(d+1)}-\lgt{d},\qquad d=1,2,\ldots,9,
	  \label{Eq:benford}
\end{equation}
where $\Pr{A}$ denotes the probability of an event $A$.
The distribution of the leading digit $d_0^{}$ depends, of course, on
that of $X$.  Although there are many continuous classical
distributions, to the author's knowledge, none of them gives the \Blaw
exactly.  The purpose of this note is to construct a continuous
distribution of $X$ on $(0,\infty)$ such that $S(X)$ follows \Blaw
exactly.

Hereafter we consider only the case $X>0$, and assume that the
base\,(radix), which will be denoted by $b$, is not necessarily 10,
but instead any integer $b\geq 2$.  Like the decimal case, we define
``scientific notation'' for the base $b$ as
\begin{equation}
	  X=b^M\cdot S(X). 	   \label{Eq:scientific_b} 
\end{equation}
In this case, the significand $S(X)$ and the exponent $M$ are given by
\[
   S(X)=b^F, \qquad  M=\lfloor \lgb{X} \rfloor, ~~~F=\{\lgb{X}\},
\] 
where the symbols $\lfloor ~~ \rfloor$ and $\{ ~~ \}$ are the integer
and the fractional parts of the argument, respectively.

For the base $b$, we extend \Blaw as 
\begin{equation}
   \Pr{1\leq S(X)\leq s}=\lgb{s},\qquad 1\leq s < b,
	  \label{Eq:s_blaw}
\end{equation}
which is often called the {\em strong Benford's law} in base $b$.  We
will simply refer to the strong law as Benford's law.  Of course, this
definition implies \eqref{Eq:benford} when $b=10$.  Since
$F=\lgb{S}$, we can easily see that the condition
\eqref{Eq:s_blaw} is equivalent to
\begin{equation}
   \Pr{0\leq F \leq \sigma}=\sigma, \qquad 0\leq \sigma < 1,
\label{Eq:uniform}
\end{equation}
where $\sigma=\lgb{s}$.  This means that the fractional part of
$\log_{\,b}X$ is uniformly distributed on $[\,0,\,1)$.  The uniformity
of $F$ is a key concept in Benford analysis.
\section{Piecewise continuous and continuous distributions}
Let $\PS{s}$ be the probability density function (in short pdf) of the
significand $S$.  To derive $\PS{s}$ we differentiate
\eqref{Eq:s_blaw} with respect to $s$.  Then we have
\begin{equation}
   P_S^{}(s)=\frac{1}{s\,\ln b}, \qquad s\in [1,\,b\,).
	  \label{Ex:pdfofs} 
\end{equation} 
If the range of $X$ is $[1,\,b\,)$ then this is also the pdf of $X$,
since in this case $X=S(X)$.  Next we extend the range of $X$ to a
wider range.

Let the range of $X$ be $[\,b^{m_0},\,b^{m_1}\,]\,(m_0<m_1)$, and
$(p_m^{})$ be the sequence satisfying
\begin{equation}
   p_m^{}\geq 0,~~~~m=m_0,m_0+1,\ldots, m_1,\qquad \sum_{m=m_0}^{m_1}
	  p_m=1. \label{Eq:defpM} 
\end{equation}
Using the sequence, we define the pdf of $X$ as
\begin{equation}
   P_X^{}(x)=\sum_{m=m_0}^{m_1}
	  \frac{p_m^{}}{x \ln b}\,\chi_m^{}(x), \qquad 0<x<\infty,
	  \label{Eq:pmpdf}
\end{equation}
where
\[
   \chi_m^{}(x)=\begin{cases}
				   1,  & ~~x\in [\,b^m,\,b^{m+1}\,),\\
				   0,  & ~~{\rm otherwise}.
				\end{cases}
\]
This function clearly satisfies the requirement for pdf's
\[
  \PX{x}\geq 0,\qquad \int_{-\infty}^{\infty}\PX{x}\,\d x=1.
\]
The random variable $X$ with the density \eqref{Eq:pmpdf} is clearly
Benford, since for some $s\,(1\leq s < b)$
\begin{equation}
   \Pr{1\leq S(X)\leq s}=\sum_{m=m_0^{}}^{m_1^{}} \int_{b^m}^{b^m
	  s}\PX{x}\,\d x=\sum_{m=m_0^{}}^{m_1^{}} \int_{b^m}^{b^m
	  s}\frac{p_m^{}}{x \ln b}\,\d x=\lgb{s}.
\end{equation}
Thus we have easily constructed a distribution of $X$ on the wider
range that satisfies \eqref{Eq:s_blaw}.

Next we consider a distribution on $(0, \infty)$.  In order to extend
the domain of $P_X^{}(x)$ to $\mathbb{R}^+$, we must set
$m_0=-\infty,\,m_1=+\infty$.  Moreover, if such a distribution to be
continuous then
\begin{equation}
   \cdots =p_{-2}^{}= p_{-1}^{}= p_0^{}=p_{1}^{}=p_2^{}=\cdots,
	  \label{Eq:scale} 
\end{equation}
which clearly contradicts \eqref{Eq:defpM}.  By the way, since
$p_m^{}$ means
\begin{equation}
   p_m^{}=\Pr{b^m \leq X <b^{m+1}},\qquad m=0,\pm 1,\pm2,\ldots,
\end{equation}
condition \eqref{Eq:scale} means
\[
   \cdots =\Pr{b^{-1}\leq X <1} = \Pr{1\leq X <b}=   \Pr{b\leq X
   <b^2}=   \cdots.
\]
This is just the hypothesis of {\em scale invariance} for a scale
factor $b$.  The hypothesis is often appeared in numerous articles on
Benford's law, since Pinkham \cite{pinkham} derived \Blaw from this
hypothesis.
However, Knuth \cite{knuth} showed that this attractive hypothesis is
not possible for all bases $b\geq 2$.  Conversely, it was shown by
Hamming \cite{hamming} that if $X$ is Benford, that is, $S(X)$ has the
distribution satisfying \eqref{Eq:s_blaw} then for any constants $c>0$
the significand $S(c X)$ also has that distribution.

We have failed to construct a continuous density on $(0,\infty)$ that
gives the law \eqref{Eq:s_blaw}.  Here we propose another way of
constructing such density.  Let us introduce the new variable
\begin{equation}
   Y=\ln X,
\end{equation}
and the pdf of $Y$ be $P_Y^{}(y)$.  From the well-known formula of probability
theory, we have
\begin{equation}
  \PY{y}=\left|\frac{\d x}{\d y}\right|\PX{x}=\e^y \PX{\e^y}
	 \mbox{~~~or~~~} \PX{x}=x^{-1}\PY{\ln x},
	 \label{Eq:transform}
\end{equation}
where $y=\ln x$ \cite{capinski}.  Using $P_Y^{}(y)$, we have for
$\sigma=\log_{\,b}s$
\begin{equation}
   \begin{split}
	  \Pr{0\leq F\leq \sigma}&=\Pr{1\leq S(X)\leq s}\\
	  &=\sum_{m=-\infty}^\infty \Pr{b^m \leq X \leq b^m s}\\
	  &=\sum_{m=-\infty}^\infty \Pr{m \ln b \leq Y \leq  (m+\sigma)\ln b }\\
	  &=\sum_{m=-\infty}^\infty \int_{m \ln b}^{(m+\sigma)\ln b} \PY{y}\,
	  \d y. \label{Eq:prsx}
   \end{split}
\end{equation}
To obtain $P_F^{}(\sigma)$, the pdf of $F$, we differentiate
\eqref{Eq:prsx} with respect to $\sigma$ under the integral sign to obtain
\begin{equation}
   \begin{split}
	  P_F(\sigma )&=\frac{\d}{\d \sigma}\, \Pr{0\leq F\leq \sigma}\\
	  &= \sum_{m=-\infty}^\infty 
	  \frac{\d}{\d \sigma}\left(  \int_{m \ln b}^{(m+\sigma)\ln b}
	  \PY{y}\,  {\d}y \right),\\
	  &= \ln b\,\sum_{m=-\infty}^\infty 
	  \PY{(m+\sigma)\ln b}. \label{Eq:pdf_F} 
   \end{split}
\end{equation}
This is just the trapezoidal approximation to the infinite integral
\begin{equation}
   \int_{-\infty}^{\infty}\PY{y}\,{\d}y=1   \label{Eq:integral_pdf_y}
\end{equation}
with the step size $h=\ln b$.  Therefore, if the trapezoidal rule
\eqref{Eq:pdf_F} gives the exact value of the integral
\eqref{Eq:integral_pdf_y}, then the distribution of the fractional
part $F=\{\log_{\,b} X\}$ is uniformly distributed on $[\,0,\,1\,)$,
and as a result the distribution of $S(X)$ is Benford.

Here we consider the class of pdf functions which satisfies
\begin{equation}
   \int_{-\infty}^{\infty}\PY{y}\,{\d}y=\ln b\,\sum_{m=-\infty}^\infty 
	  \PY{(m+\sigma)\ln b}.  \label{Eq:intequaltrap}
\end{equation}
\section{Trapezoidal rule and Benford distribution}
Let us consider the integral
\begin{equation}
   I=\int_a^b f(x)\,\d x.  \label{Eq:integral}
\end{equation}
To approximate the integral numerically many formulas were developed.
Among the formulas, the trapezoidal rule is the most elementary
formula.  As is well known, the formula cannot be expected to give an
accurate result for an integral over a large\,(finite) interval; this
formula gives the exact result only for the case that the integrand
$f(x)$ is a linear or piecewise linear function.  Therefore, the
formula is seldom used in practice \cite{evans}.  It is shown,
however, that if the interval is infinite or semi-infinite and $f(x)$
is an analytic function on $\mathbb{R}$ or $\mathbb{R}^+$, then the
convergence is tremendously fast \cite{trefethen}.  Moreover, Sugihara
\cite{sugihara} proved that there exists a class of functions for
which the trapezoidal rule gives the exact value of $I$, when
$a=-\infty,\,b=+\infty$.  Here we show the result by Sugihara:\\[2mm]

\noindent {\bf Theorem\,(Sugihara \cite{sugihara})} Let $f(z)$ be a
function that satisfies the following three conditions:
\begin{equation}
   \begin{split}
	  &\lefteqn{{\rm (a)}~\left|\int_{-\infty}^\infty f(x)\,\d x
	  \right|<\infty.}\\ 
	  &{\rm (b)}~f(z)~\mbox{is entire, that is, holomorphic for
	  all}~z\in{\mathbb C}.  \hspace{5.5cm}\\ 
	  &{\rm (c)}~f(z)~\mbox{is exponential type}~A.
   \end{split}  \label{Eq:trapequalint0}
\end{equation}  
Then for all $0<h<2\pi/A$
\begin{equation}
   \int_{-\infty}^\infty f(x)\,\d x =h \sum_{m=-\infty}^\infty f(mh).
	  \label{Eq:trapeqint}
\end{equation}
~\\ As an example of the functions belonging to this class with
$A=1$, we show $f(z)=\sinc{z}(=\sin z/z)$.  For this function we have
\begin{equation}
   h \sum_{m=-\infty}^\infty \sinc{m h}=\int_{-\infty}^\infty
	 \sinc{x}\,\d x\,(= \pi)
\end{equation}
for all $0< h <2\pi$.

In the present situation, since the integrand is a pdf, we must modify
the above conditions as follows:
\begin{equation}
\begin{split}
   & \lefteqn{{\rm (a)}~~ \int_{-\infty}^\infty f(x)\,\d x =1}\\
   & {\rm (b)}~~f(z)~\mbox{is entire},\mbox{~and~}
   f(x)\geq 0~\mbox{for}~x\in {\mathbb R}.\hspace{7.0cm}\\
   & {\rm (c)}~~ f(z)~\mbox{is exponential type}~A.
\end{split}\label{Eq:trapequalint1}
\end{equation}

We now return to the topic of Benford's law.  From now on the
integrand is not $f(z)$ but $\PY{y}$, and the argument is not $x$ but
$y$.  As an example of the functions that satisfy the three conditions
in \eqref{Eq:trapequalint1}, we can show the following function:
\begin{equation}
   \PY{y}=\left(\frac{a}{\pi}\right){\rm sinc}^2 (a\,y),\qquad a>0. 
	  \label{Eq:PYsinc}
\end{equation}
It can easily be shown that this function satisfies the conditions (a)
and (b).  It is also shown that the function is exponential type
$2a$\,\cite{rahman}, since on the imaginary axis
\[
   \lim_{v\to \pm\infty}\frac{\log \left| \PY{{\rm
   i}\,v}\right|}{|v|}=2a.
\]
Therefore, for the $\PY{y}$ given by \eqref{Eq:PYsinc} and for the $b$
in the range
\begin{equation}
   0<\ln b < \pi/a,   \label{Eq:h_range}
\end{equation}
eq.\eqref{Eq:intequaltrap} is valid.  Thus the random variable $X$
with the density
\begin{equation}
   \PX{x}=x^{-1}\PY{\ln x}=\left(\frac{a}{\pi x}\right){\rm sinc}^2 (a
	  \ln x), \qquad  x>0  \label{Eq:x_density}
\end{equation}
is Benford, that is, $S(X)$ satisfy the condition \eqref{Eq:s_blaw}
for all $b$ satisfying
\begin{equation}
   2\leq b<\e^{\pi/a}.  \label{Eq:range_of_base}
\end{equation}
In particular, in order that the variable $X$ is Benford for $b=10$,
we must keep $a$ in the range
\begin{equation}
   0<a<\frac{\pi}{\ln 10}=1.364\cdots.
\end{equation}

Here we consider the mean and variance of the random variable $X$ with
the distribution \eqref{Eq:x_density}.  Let $\lambda\geq 1$ be an
integer, then $\lambda$\,th moment is
\begin{equation}
   \mathbb{E}\,[X^\lambda]=\int_0^\infty x^\lambda P_X^{} (x)\,{\rm d}x
   = \frac{a}{\pi}\int_{-\infty}^{\infty} {\rm sinc}^2 (a
   y)\,\e^{\lambda y}\, \d y,
\end{equation}
where $\mathbb{E}\,[~\cdot~]$ denotes the expected value of the argument.
Since in this equation the integrand diverges as $y\to +\infty$, then
the mean and variance of $X$ do not exist.
\section{Conclusion}
We have developed a continuous distribution on $(0,\infty)$ of the
random variable whose significand obeys \Blaw exactly by using
Sughihara's theory on the trapezoidal rule \cite{sugihara}.  It has
been shown that the random variable with the distribution does not
have the mean and the variance.
\end{document}